\documentstyle[12pt,leqno]{article}
\newtheorem{theorem}{Theorem}
\newtheorem{definition}[theorem]{Definition}

\newtheorem{example}[theorem]{Example}

\catcode`\@=11

\def\ps@myheadings{\let\@mkboth\@gobbletwo
\def\@oddhead{\ifnum\count0=1 \hfill\else
\rightmark \hfil \rm\thepage\fi}%
\def\@oddfoot{\ifnum\count0=1 \hfill \rm 1 \hfill \else
\hfill\fi}
\def\@evenhead%
{\rm\leftmark\hfil\rm\thepage}%
\def\@evenfoot{}\def\sectionmark##1{}
\def\subsectionmark##1{}}

\def\@begintheorem#1#2{\it \trivlist \item[\hskip
 \labelsep{\bf #1\ #2.}]}
\def\@opargbegintheorem#1#2#3{\it \trivlist\item[\hskip%
 \labelsep{\bf #1\ #2.\ (#3)}]}
\def\@endtheorem{\endtrivlist}

\def\@listI{\leftmargin\leftmargini \parsep 1pt plus 2.5pt
 minus 1pt\topsep 5pt plus 2pt minus 3pt%
 \itemsep 0pt plus 2.5pt minus 1pt}
\let\@listi\@listI
\@listi

\def\@sect#1#2#3#4#5#6[#7]#8{\ifnum #2>\c@secnumdepth%
 \def \@svsec {}\else \refstepcounter {#1}\edef \@svsec%
 {\csname the#1\endcsname. \hskip .1em }\fi \@tempskipa%
 #5\relax \ifdim \@tempskipa >\z@ \begingroup #6\relax%
 \@hangfrom {\hskip #3\relax \@svsec }{\interlinepenalty%
 \@M #8.\par }\endgroup \csname #1mark\endcsname {#7}%
 \addcontentsline {toc}{#1}{\ifnum #2>\c@secnumdepth%
 \else \protect \numberline {\csname the#1\endcsname. }%
 \fi #7}\else \def \@svsechd {#6\hskip #3\@svsec #8.%
 \csname #1mark\endcsname {#7}\addcontentsline {toc}{#1}%
 {\ifnum #2>\c@secnumdepth \else \protect \numberline%
 {\csname the#1\endcsname. }\fi #7}}\fi \@xsect {#5}}

\def\section{\@startsection {section}{1}{\z@ }%
 {-3.5ex plus -1ex minus -.2ex}{2.3ex plus .2ex}{\bf }}

\def\thebibliography#1{%
 \section *{References.\@mkboth {REFERENCES}{REFERENCES}}%
 \list {[\arabic {enumi}]}{\settowidth \labelwidth {[#1]}%
 \leftmargin \labelwidth \advance \leftmargin \labelsep %
 \usecounter {enumi}} \def \newblock %
 {\hskip .11em plus .33em minus -.07em} \sloppy \clubpenalty 4000%
 \widowpenalty 4000 \sfcode`\.=1000\relax}

\def\@maketitle{%
 \newpage \null \vskip 2em
 \begin{center}
{\Large\bf \@title \par }
 \vskip 1.5em
 {\large \lineskip .5em
 \begin {tabular}[t]{c}\@author
 \end{tabular}\par}
 \end{center}
  \vskip .8em}

\def\abstract{%
\if@twocolumn \section *{Abstract}
 \else \small\quotation\noindent{\bf Abstract.}\fi}

\catcode`\@=13

\font\tenbold=msbm10 scaled \magstep1
\font\sevenbold=msbm7 scaled \magstep1
\font\fivebold=msbm5 scaled \magstep1
\newfam\boldfam 
\textfont\boldfam=\tenbold
\scriptfont\boldfam=\sevenbold
\scriptscriptfont\boldfam=\fivebold

\def\ot{\otimes}

\def\prada#1#2{#1 + \cdots + #2}
\def\Rada#1#2#3{#1_{#2},\dots,#1_{#3}}
\def\pa{\partial}
\def\znamenko#1{{(-1)^{#1}\cdot}}
\def\doubless#1#2{{
\def\arraystretch{.5}
\begin{array}{c}
\mbox{\scriptsize $\scriptstyle #1$}
\\
\mbox{\scriptsize $\scriptstyle #2$}
\end{array}\def\arraystretch{1}
}}

\def\Ass{{\it Ass}}
\def\Lie{{\it Lie}}
\def\hpalg{{h{\cal P}}}
\def\ainfty{{\rm A}(\infty)}
\def\cinfty{{\rm C}(\infty)}
\def\linfty{{\rm L}(\infty)}
\def\P{{\cal P}}
\def\F{{\cal F}}
\def\M{{\cal M}}

\def\oper{{
\unitlength=.500000pt
\begin{picture}(25.00,25.00)(30.00,10.00)
\thicklines
\put(40.00,20.00){\makebox(0.00,0.00){$\bullet$}}
\put(40.00,20.00){\line(1,-1){20.00}}
\put(40.00,20.00){\line(-1,-1){20.00}}
\put(40.00,40.00){\line(0,-1){20.00}}
\end{picture}
}}

\def\opern{{
{
\unitlength=1.000000pt
\thicklines
\begin{picture}(40.00,28.00)(0.00,20.00)
\put(20.00,0.00){\makebox(0.00,0.00)[l]{$\cdots$}}
\put(20.00,20.00){\makebox(0.00,0.00){$\bullet$}}
\put(20.00,20.00){\line(1,-1){20.00}}
\put(20.00,20.00){\line(-1,-2){10.00}}
\put(20.00,20.00){\line(-1,-1){20.00}}
\put(20.00,40.00){\line(0,-1){20.00}}
\end{picture}}
}}

\def\opernm{{
{
\unitlength=.6pt
\thicklines
\begin{picture}(40.00,30.00)(0.00,-3.00)
\put(13.00,0.00){\makebox(0.00,0.00)[l]{$\cdots$}}
\put(20.00,20.00){\makebox(0.00,0.00){$\bullet$}}
\put(20.00,20.00){\line(1,-1){20.00}}
\put(20.00,20.00){\line(-1,-2){10.00}}
\put(20.00,20.00){\line(-1,-1){20.00}}
\put(20.00,40.00){\line(0,-1){20.00}}
\end{picture}}
}}

\def\slash{
{
\thicklines
\unitlength=1.20000pt
\begin{picture}(10.00,20.00)(0.00,7.00)
\put(0.00,0.00){\line(1,2){10.00}}
\end{picture}}
}

\title{Homotopy Algebras via Resolutions of Operads}

\author{Martin Markl%
\thanks{The author was supported by the
grant GA AV \v CR 1019804. This paper is in final form and no version
of it will be submitted for publication elsewhere.}}

\begin{document}
\bibliographystyle{alpha}
\baselineskip 18pt plus 2pt minus 1 pt

\maketitle

All algebraic objects in this note will be considered over a fixed
field $\bf k$ of characteristic zero. If not stated otherwise, all
operads live in the category of differential graded vector spaces over
${\bf k}$. For standard terminology concerning operads, algebras over
operads, etc., see either the original paper by
J.P.~May~\cite{may:1972}, or an overview~\cite{loday:bourbaki}.

The aim of this note is mainly to advocate our approach to homotopy
algebras based on the minimal model of an operad.
Our intention is
to expand it to a paper on homotopy properties of
the category of homotopy algebras; some possible results in this
direction are indicated in Section~\ref{lights}.

I would like to express my thanks to Rainer Vogt for his kind
invitation to Osnabr\"uck. I also owe my thanks to Jim Stasheff for
careful reading the manuscript and many useful comments.

\section{Motivations}  

\begin{example}{\rm
\label{the}
An associative algebra consists of a vector space $A$ together with a
multiplication $\mu :A \ot A \to A$, which is supposed to be
associative:
\[
\mu(\mu(a,b),c) = \mu(a,\mu(b,c));
\]
we neglect the role of a unit.

We can `visualize' $\mu$ as an `operation' with two inputs and one output,
$\mu =\oper$\ .
The associativity can be then depicted as
\begin{center}
{
\unitlength=0.5pt
\thicklines
\begin{picture}(190.00,60.00)(0.00,0.00)
\put(90.00,40.00){\makebox(0.00,0.00)[l]{$=$}}
\put(170.00,20.00){\makebox(0.00,0.00){$\bullet$}}
\put(150.00,40.00){\makebox(0.00,0.00){$\bullet$}}
\put(20.00,20.00){\makebox(0.00,0.00){$\bullet$}}
\put(40.00,40.00){\makebox(0.00,0.00){$\bullet$}}
\put(150.00,60.00){\line(0,-1){20.00}}
\put(170.00,20.00){\line(-1,-1){20.00}}
\put(150.00,40.00){\line(1,-1){40.00}}
\put(110.00,0.00){\line(1,1){40.00}}
\put(40.00,40.00){\line(1,-1){40.00}}
\put(20.00,20.00){\line(1,-1){20.00}}
\put(40.00,40.00){\line(-1,-1){40.00}}
\put(40.00,60.00){\line(0,-1){20.00}}
\end{picture}}
\end{center}
This means that the operad $\Ass$ describing associative algebras in
the category of differential graded vector spaces can be
presented as the quotient
\[
{
{\cal F}\left(\ \ \oper\ \right)\slash
}{\left(\hskip 2mm
{
\unitlength=0.5pt
\thicklines
\begin{picture}(190.00,40.00)(0.00,25.00)
\put(90.00,40.00){\makebox(0.00,0.00)[l]{$-$}}
\put(170.00,20.00){\makebox(0.00,0.00){$\bullet$}}
\put(150.00,40.00){\makebox(0.00,0.00){$\bullet$}}
\put(20.00,20.00){\makebox(0.00,0.00){$\bullet$}}
\put(40.00,40.00){\makebox(0.00,0.00){$\bullet$}}
\put(150.00,60.00){\line(0,-1){20.00}}
\put(170.00,20.00){\line(-1,-1){20.00}}
\put(150.00,40.00){\line(1,-1){40.00}}
\put(110.00,0.00){\line(1,1){40.00}}
\put(40.00,40.00){\line(1,-1){40.00}}
\put(20.00,20.00){\line(1,-1){20.00}}
\put(40.00,40.00){\line(-1,-1){40.00}}
\put(40.00,60.00){\line(0,-1){20.00}}
\end{picture}}
\hskip 2mm
\right)}
\]
where ${\cal F}\left(\ \ \oper\ \right)$
is the free non-$\Sigma$ operad on the operation \oper\ representing
the product $\mu$ and
\[
{\left(\hskip 2mm
{
\unitlength=0.5pt
\thicklines
\begin{picture}(190.00,40.00)(0.00,25.00)
\put(90.00,40.00){\makebox(0.00,0.00)[l]{$-$}}
\put(170.00,20.00){\makebox(0.00,0.00){$\bullet$}}
\put(150.00,40.00){\makebox(0.00,0.00){$\bullet$}}
\put(20.00,20.00){\makebox(0.00,0.00){$\bullet$}}
\put(40.00,40.00){\makebox(0.00,0.00){$\bullet$}}
\put(150.00,60.00){\line(0,-1){20.00}}
\put(170.00,20.00){\line(-1,-1){20.00}}
\put(150.00,40.00){\line(1,-1){40.00}}
\put(110.00,0.00){\line(1,1){40.00}}
\put(40.00,40.00){\line(1,-1){40.00}}
\put(20.00,20.00){\line(1,-1){20.00}}
\put(40.00,40.00){\line(-1,-1){40.00}}
\put(40.00,60.00){\line(0,-1){20.00}}
\end{picture}}
\hskip 2mm
\right)}
\]
is the operadic ideal generated by the element
\[
{
{
\unitlength=0.5pt
\thicklines
\begin{picture}(190.00,60.00)
\put(90.00,40.00){\makebox(0.00,0.00)[l]{$-$}}
\put(170.00,20.00){\makebox(0.00,0.00){$\bullet$}}
\put(150.00,40.00){\makebox(0.00,0.00){$\bullet$}}
\put(20.00,20.00){\makebox(0.00,0.00){$\bullet$}}
\put(40.00,40.00){\makebox(0.00,0.00){$\bullet$}}
\put(150.00,60.00){\line(0,-1){20.00}}
\put(170.00,20.00){\line(-1,-1){20.00}}
\put(150.00,40.00){\line(1,-1){40.00}}
\put(110.00,0.00){\line(1,1){40.00}}
\put(40.00,40.00){\line(1,-1){40.00}}
\put(20.00,20.00){\line(1,-1){20.00}}
\put(40.00,40.00){\line(-1,-1){40.00}}
\put(40.00,60.00){\line(0,-1){20.00}}
\end{picture}}
}
\]
corresponding to the associativity axiom. 
As differential graded vector spaces, the pieces $\Ass(n)$, $n\geq 1$,
of the operad $\Ass$ are
concentrated in degree zero and have trivial differentials.

Homotopy versions of associative algebras are so-called
$\ainfty$-algebras (also called strongly homotopy associative
algebras), 
introduced by 
J.~Stasheff~\cite{stasheff:TAMS63}. An
$\ainfty$-algebra $A = (A,\pa,\mu_2,\mu_3,\ldots)$ consists of a
differential graded 
vector space $A = (A,\pa)$ 
and multilinear operations $\mu_n : A^{\ot n} \to A$
of degree $n-2$ that satisfy the following infinite set of axioms:
\begin{eqnarray}
\label{pentagon}
\mu_2(\mu_2(a,b)c) - \mu_2(a,\mu_2(b,c)) \hskip -2mm &=& \hskip -3mm [\mu_3,\pa](a,b,c),
\\
\nonumber 
\mu_3(\mu_2(a,b),c,d) -
\mu_3(a,\mu_2(b,c),d)+\mu_3(a,b,\mu_2(c,d))- \hskip1cm &&
\\ 
\nonumber 
- \mu_2(\mu_3(a,b,c),d) -\znamenko{|a|}\mu_2(a,\mu_3(b,c,d))
\hskip -2mm &=& \hskip -3mm
[\mu_4,\pa](a,b,c,d),
\\
\nonumber 
&\vdots&
\\
\nonumber 
\sum_{\doubless{i+j=n+1}{i,j \geq 2}} \hskip -2mm
\sum_{s=0}^{n-j}
\znamenko\epsilon
\mu_i(\Rada a1s,\mu_j(\Rada a{s+1}{s+j}),\Rada a{s+j+1}n)
\hskip -2mm &=& \hskip -3mm
[\mu_n,\pa](\Rada a1n),
\end{eqnarray}
where $[\mu_n,\pa](\Rada a1n)$ denotes the induced differential
in the endomorphism complex ${\rm Hom}_{\bf k}(A^{\ot n},A)$, 
\begin{eqnarray*}
\lefteqn{
[\mu_n,\pa](\Rada a1n) :=}
\\&&:= 
\sum_{1\leq s \leq n} 
\znamenko{\prada{|a_1|}{|a_{s-1}|}}\mu_n(a_1,\ldots,\partial
a_s,\ldots,a_n) -\znamenko n\partial \mu_n(\Rada a1n).
\end{eqnarray*}
The sign is given by
\[
\epsilon :=j+s(j+1)+j(\prada{|a_1|}{|a_{s-1}|}).
\]
The first axiom of~(\ref{pentagon}) 
expresses the homotopy associativity of the product
$\mu_2$. 
The second axiom is
a linearized version of the Mac~Lane pentagon condition, with the terms in
the left hand side corresponding to the edges of the pentagon:

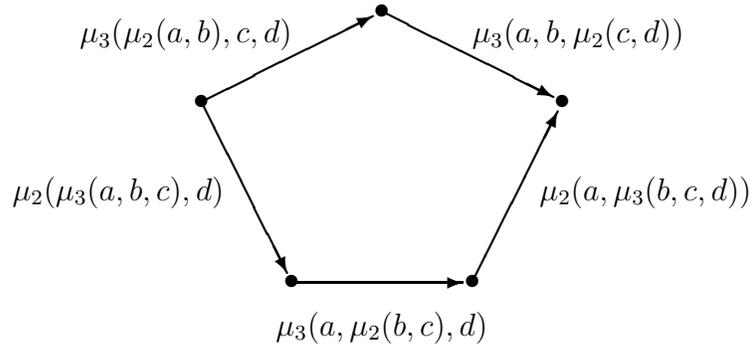
\begin{center}
\setlength{\unitlength}{1.2cm}
\thicklines
\begin{picture}(4,4)(-2,-0.5)

\put(-1,0){\vector(1,0){1.9}}
\put(1,0){\vector(1,2){.95}}
\put(0,3){\vector(2,-1){1.9}}
\put(-2,2){\vector(2,1){1.9}}
\put(-2,2){\vector(1,-2){.95}}

\put(1,0){\makebox(0,0){$\bullet$}}
\put(2,2){\makebox(0,0){$\bullet$}}
\put(0,3){\makebox(0,0){$\bullet$}}
\put(-2,2){\makebox(0,0){$\bullet$}}
\put(-1,0){\makebox(0,0){$\bullet$}}

\put(0,-.5){\makebox(0,0){{$\mu_3(a,\mu_2(b,c),d)$}}}
\put(-1.75,1){\makebox(0,0)[r]{{$\mu_2(\mu_3(a,b,c),d)$}}}
\put(1.75,1){\makebox(0,0)[l]{{$\mu_2(a,\mu_3(b,c,d))$}}}
\put(-1,2.75){\makebox(0,0)[r]{{$\mu_3(\mu_2(a,b),c,d)$}}}
\put(1,2.75){\makebox(0,0)[l]{{$\mu_3(a,b,\mu_2(c,d))$}}}

\end{picture}
\end{center}

This miraculous phenomenon is explained by 
the fact that the operad $\ainfty$ is isomorphic to
the operad $\{CC_*(K_n)\}_{n\geq 1}$ of cellular chains
of a certain cellular operad 
$\{K_n\}_{n\geq 1}$,
introduced by J.~Stasheff in his pioneering 
work~\cite{stasheff:TAMS63}. The spaces $K_n$ are now called
associahedra; the 
polyhedron $K_4$ that parameterizes 4-ary operations  
is the pentagon mentioned above.

A more straightforward description of $\ainfty$ can be
obtained from the axioms as follows. The differential graded 
operad $\ainfty$
is free as an operad,
\[
\ainfty := {\cal F}\left(
{
\unitlength=.500000pt
\thicklines
\begin{picture}(180.00,40.00)(0.00,15.00)
\put(180.00,20.00){\makebox(0.00,0.00){$\cdots$}}
\put(110.00,20.00){\makebox(0.00,0.00){,}}
\put(50.00,20.00){\makebox(0.00,0.00){,}}
\put(20.00,20.00){\makebox(0.00,0.00){$\bullet$}}
\put(140.00,20.00){\makebox(0.00,0.00){$\bullet$}}
\put(80.00,20.00){\makebox(0.00,0.00){$\bullet$}}
\put(140.00,20.00){\line(0,1){0.00}}
\put(140.00,40.00){\line(0,-1){20.00}}
\put(140.00,20.00){\line(1,-2){10.00}}
\put(140.00,20.00){\line(-1,-2){10.00}}
\put(140.00,20.00){\line(1,-1){20.00}}
\put(120.00,0.00){\line(1,1){20.00}}
\put(80.00,40.00){\line(0,-1){40.00}}
\put(80.00,20.00){\line(1,-1){20.00}}
\put(60.00,0.00){\line(1,1){20.00}}
\put(20.00,20.00){\line(1,-1){20.00}}
\put(20.00,20.00){\line(-1,-1){20.00}}
\put(20.00,20.00){\line(0,1){0.00}}
\put(20.00,40.00){\line(0,-1){20.00}}
\end{picture}}
\hskip 3mm \right),\
\deg \left(
\underbrace{\opernm}_{\mbox{$n$-times}}
\right)
= n-2
\]
with the differential given on generators by
\begin{center}
\unitlength 1.50mm
\thicklines
\begin{picture}(44.50,20.00)(0,15)
\put(22.00,25.00){\makebox(0,0)[cc]{$
\pa\left(\
\opern\
\right)
 := \displaystyle%
\sum_{\doubless{i+j=n+1}{i,j \geq 2}}
\sum_{s=0}^{n-j}
\znamenko{j+s(j+1)}
\left(\rule{0mm}{17mm}\hskip40mm \right)$}}
\put(5,0){\put(35.17,35.00){\line(0,-1){5.00}}
\put(35.17,30.00){\line(-6,-5){6.00}}
\put(35.17,30.00){\line(-1,-2){2.50}}
\put(35.17,29.67){\line(1,-3){3.11}}
\put(38.28,20.33){\line(-1,-1){6.11}}
\put(38.33,20.00){\line(-2,-5){2.27}}
\put(38.33,20.00){\line(1,-1){5.50}}
\put(35.17,29.83){\line(2,-1){9.33}}
\put(39.17,15.83){\makebox(0,0)[cc]{$\cdots$}}
\put(39.50,25.33){\makebox(0,0)[cc]{$\cdots$}}
\put(34.67,25.33){\makebox(0,0)[cc]{$\cdots$}}
\put(39.17,22.00){\makebox(0,0)[lc]{$s+1$-th input}}
\put(35.17,29.8){\makebox(0,0)[cc]{$\bullet$}}
\put(38.33,20.00){\makebox(0,0)[cc]{$\bullet$}}
}
\end{picture}
\end{center}
and extended by the derivation property.
$\ainfty$-algebras occur in Nature as chain algebras of loop spaces,
a result which is today classical, 
see again~\cite{stasheff:TAMS63}. 
}
\end{example}

\begin{example}{\rm
A Lie algebra is a vector space $L$ with an antisymmetric
product $[-,-]: L\ot L \to L$ (the `bracket'), satisfying the Jacoby
identity:
\[
[a,[b,c]] + [b,[c,a]] + [c,[a,b]] = 0.
\]
Thus the operad $\Lie$ for Lie algebras is the quotient
\[
{\cal F}\left(\ \ \oper\ \right) \slash
\left(
{
\unitlength=1.00000pt
\thicklines
\begin{picture}(170.00,25.00)(-5.00,25.00)
\put(110.00,30.00){\makebox(0.00,0.00){$+$}}
\put(50.00,30.00){\makebox(0.00,0.00){$+$}}
\put(120.00,0.00){\makebox(0.00,0.00){$3$}}
\put(80.00,0.00){\makebox(0.00,0.00){$3$}}
\put(40.00,0.00){\makebox(0.00,0.00){$3$}}
\put(160.00,0.00){\makebox(0.00,0.00){$2$}}
\put(60.00,0.00){\makebox(0.00,0.00){$2$}}
\put(20.00,0.00){\makebox(0.00,0.00){$2$}}
\put(140.00,0.00){\makebox(0.00,0.00){$1$}}
\put(100.00,0.00){\makebox(0.00,0.00){$1$}}
\put(0.00,0.00){\makebox(0.00,0.00){$1$}}
\put(30.00,20.00){\makebox(0.00,0.00){$\bullet$}}
\put(90.00,20.00){\makebox(0.00,0.00){$\bullet$}}
\put(150.00,20.00){\makebox(0.00,0.00){$\bullet$}}
\put(140.00,30.00){\makebox(0.00,0.00){$\bullet$}}
\put(80.00,30.00){\makebox(0.00,0.00){$\bullet$}}
\put(20.00,30.00){\makebox(0.00,0.00){$\bullet$}}
\put(150.00,20.00){\line(0,1){0.00}}
\put(140.00,10.00){\line(1,1){10.00}}
\put(140.00,50.00){\line(0,-1){20.00}}
\put(140.00,30.00){\line(1,-1){20.00}}
\put(120.00,10.00){\line(1,1){20.00}}
\put(80.00,50.00){\line(0,-1){20.00}}
\put(90.00,20.00){\line(-1,-1){10.00}}
\put(80.00,30.00){\line(1,-1){20.00}}
\put(60.00,10.00){\line(1,1){20.00}}
\put(20.00,10.00){\line(1,1){10.00}}
\put(20.00,30.00){\line(1,-1){20.00}}
\put(0.00,10.00){\line(0,1){0.00}}
\put(20.00,30.00){\line(-1,-1){20.00}}
\put(20.00,50.00){\line(0,-1){20.00}}
\end{picture}}
\right)
\]
Here ${\cal F}\left(\ \oper \ \right)$ now denotes the free {\em
symmetric\/} ($\Sigma$-) operad on one antisymmetric bilinear
operation; the labels of the inputs encode the action of the symmetric
group in a standard manner.

Homotopy versions of Lie algebras are strongly homotopy  Lie algebras
(also called sh Lie algebras or $\linfty$-algebras). 
They were introduced and systematically studied
in~\cite{lada-stasheff:IJTP93}, though they had existed in the
literature, in 
various disguises, even before; see also~\cite{lada-markl:CommAlg95}. 

An $\linfty$
algebra is a differential graded vector space $L = (L,\pa)$, together
with a set $\{l_n\}_{n\geq 2}$ of graded antisymmetric operations $l_n
: L^{\ot n} \to L$ of degree $n-2$ such that the following infinite
set of axioms is satisfied for any $n\geq 2$:
\begin{eqnarray*}
\sum_%
\doubless{i+j=n+1}{i,j\geq 2}\sum_\sigma
\chi(\sigma)(-1)^{i(j-1)}l_j(l_i(\Rada a{\sigma(1)}{\sigma(i)}),\Rada
a{\sigma(i+1)}{\sigma (n)}) =&&
\\
&&\hskip -3cm =
\znamenko{n}[l_n,\pa](\Rada a1n)
\end{eqnarray*}
The summation is taken over all $(i,n-i)$-unshuffles 
\[
\sigma \in
\Sigma_n,\ 
\sigma(1) < \cdots < \sigma(i),\ 
\sigma(i+1)< \cdots < \sigma(n),
\] 
with
$n-1\geq i\geq 1$, and $\chi(\sigma)$ is a certain sign which we will
not specify here, see~\cite{lada-markl:CommAlg95}.
We write the first two axioms explicitly (though, to save
paper, 
without signs). The first axiom says that the `bracket'
$l_2$ satisfies the Jacobi
identity up to a homotopy:
\[
l_2(l_2(a,b),c) + l_2(l_2(b,c),a)+ l_2(l_2(c,a),b) = [l_3,\pa](a,b,c).
\]
The second condition reads
\def\pomoc#1#2#3#4{l_2(l_3(#1,#2,#3),#4)}
\def\Pomoc#1#2#3#4{l_3(l_2(#1,#2),#3,#4)}
\begin{eqnarray*}
\pomoc abcd + \pomoc acdb +\pomoc abdc + \pomoc bcda+&&
\\  
+\Pomoc abcd +
\Pomoc acbd + +  \Pomoc adbc+&&
\\
 + \Pomoc bcad + \Pomoc cdab 
+ \Pomoc bdac
=&& \hskip -5mm [\pa,l_4](a,b,c,d).
\end{eqnarray*}

To indicate an interpretation of the last equality, recall that
there is a Lie-analog of $K_4$, which we introduced
in~\cite{markl-shnider:ch-q-alg}. 
We denoted it by $L_4$ and called the Lie-hedron. It is
not a polyhedron, but just a graph, 
and the 10 terms in the left hand side
of the above equation correspond to the vertices (not edges!) 
of $L_4$, which is the
Peterson graph:

\begin{center}
\unitlength 2.00mm
\thicklines
{
\unitlength=7.000000pt
\begin{picture}(54.62,42.50)(0.00,0.00)
\put(23.00,9.13){\makebox(0.00,0.00)[tl]{$\Pomoc adbc$}}
\put(28.62,29.50){\makebox(0.00,0.00)[l]{$\pomoc abdc$}}
\put(35,23.50){\makebox(0.00,0.00)[bl]{$\pomoc acdb$}}
\put(33.00,12.5){\makebox(0.00,0.00)[tl]{$\Pomoc bdac$}}
\put(18.88,21.00){\makebox(0.00,0.00)[tr]{$\Pomoc acbd$}}
\put(2.00,22.88){\makebox(0.00,0.00){$\pomoc abcd$}}
\put(15.12,0.38){\makebox(0.00,0.00){$\Pomoc bcad$}}
\put(40.38,0.00){\makebox(0.00,0.00){$\pomoc bcda$}}
\put(53.62,22.62){\makebox(0.00,0.00){$\Pomoc cdab$}}
\put(27.88,40.50){\makebox(0.00,0.00){$\Pomoc abcd$}}
\put(17.88,2.50){\line(1,0){20.00}}
\put(37.88,2.50){\line(1,2){10.00}}
\put(47.88,22.50){\line(-4,3){20.00}}
\put(27.88,37.50){\line(-4,-3){20.00}}
\put(7.88,22.50){\line(1,-2){10.00}}
\put(7.88,22.50){\line(1,0){40.00}}
\put(17.88,2.50){\line(2,3){5.00}}
\put(22.88,10.00){\line(1,4){5.13}}
\put(27.88,30.50){\line(1,-4){5.13}}
\put(33.12,10.00){\line(2,-3){5.00}}
\put(22.88,10.00){\line(6,5){15.00}}
\put(33.04,10.17){\line(-5,4){15.33}}
\put(27.88,37.50){\line(0,-1){7.00}}
\put(27.88,37.50){\makebox(0.00,0.00){$\bullet$}}
\put(7.98,22.50){\makebox(0.00,0.00){$\bullet$}}
\put(17.88,2.50){\makebox(0.00,0.00){$\bullet$}}
\put(38.08,2.50){\makebox(0.00,0.00){$\bullet$}}
\put(47.71,22.50){\makebox(0.00,0.00){$\bullet$}}
\put(27.97,30.50){\makebox(0.00,0.00){$\bullet$}}
\put(17.88,22.50){\makebox(0.00,0.00){$\bullet$}}
\put(37.88,22.50){\makebox(0.00,0.00){$\bullet$}}
\put(32.97,10.10){\makebox(0.00,0.00){$\bullet$}}
\put(22.88,10.00){\makebox(0.00,0.00){$\bullet$}}
\end{picture}}

\end{center}

\hphantom{c} 
\vglue2mm
As in Example~\ref{the}, the operad $\linfty$ for strongly
homotopy Lie algebras is a differential graded operad, which is free
as an operad,
\[
\linfty := {\cal F}\left(
{
\unitlength=.500000pt
\thicklines
\begin{picture}(180.00,40.00)(0.00,15.00)
\put(180.00,20.00){\makebox(0.00,0.00){$\cdots$}}
\put(110.00,20.00){\makebox(0.00,0.00){,}}
\put(50.00,20.00){\makebox(0.00,0.00){,}}
\put(20.00,20.00){\makebox(0.00,0.00){$\bullet$}}
\put(140.00,20.00){\makebox(0.00,0.00){$\bullet$}}
\put(80.00,20.00){\makebox(0.00,0.00){$\bullet$}}
\put(140.00,20.00){\line(0,1){0.00}}
\put(140.00,40.00){\line(0,-1){20.00}}
\put(140.00,20.00){\line(1,-2){10.00}}
\put(140.00,20.00){\line(-1,-2){10.00}}
\put(140.00,20.00){\line(1,-1){20.00}}
\put(120.00,0.00){\line(1,1){20.00}}
\put(80.00,40.00){\line(0,-1){40.00}}
\put(80.00,20.00){\line(1,-1){20.00}}
\put(60.00,0.00){\line(1,1){20.00}}
\put(20.00,20.00){\line(1,-1){20.00}}
\put(20.00,20.00){\line(-1,-1){20.00}}
\put(20.00,20.00){\line(0,1){0.00}}
\put(20.00,40.00){\line(0,-1){20.00}}
\end{picture}}
\hskip 3mm \right),\
\underbrace{\opernm}_{\mbox{$n$-times}} 
\mbox{antisymmetric of degree $n-2$,}
\]
where ${\cal F}(-)$ is the free {\em symmetric\/}
($\Sigma$-)
operad. The action of the differential is
given by
\begin{center}
\setlength{\unitlength}{0.006in}%
\begin{picture}(114,190)(-110,540)
\thicklines
\put(-44,645){\makebox(0,0)[cc]{$\znamenko{n}
\pa \left(\ \opern\ \right) := \hskip -5mm\displaystyle
\sum_\doubless{i+j=n+1}{i,j\geq 2}\sum_\sigma
\chi(\sigma)(-1)^{i(j-1)}\cdot%
\left(\rule{0mm}{16mm}\hskip50mm \right)$}}
\put(60,00){
\put(160,730){\line( 0,-1){ 60}}
\put(160,670){\line(-5,-6){ 50}}
\put(110,610){\line(-1,-1){ 50}}
\put(110,610){\line(-2,-5){ 20}}
\put(110,610){\line( 4,-5){ 40}}
\put(160,670){\line( 1,-4){ 25}}
\put(185,570){\line( 0,-1){  5}}
\put(160,670){\line( 3,-2){150}}
\put(160,665){\makebox(0,0)[b]{$\bullet$}}
\put(110,607){\makebox(0,0)[b]{$\bullet$}}
\put(115,570){\makebox(0,0)[b]{$\cdots$}}
\put(235,580){\makebox(0,0)[b]{$\cdots$}}
\put( 53,540){\makebox(0,0)[b]{\scriptsize $\sigma(1)$}}
\put( 92,540){\makebox(0,0)[b]{\scriptsize $\sigma(2)$}}
\put(153,540){\makebox(0,0)[b]{\scriptsize $\sigma(i)$}}
\put(204,550){\makebox(0,0)[b]{\scriptsize $\sigma(i+1)$}}
\put(320,550){\makebox(0,0)[b]{\scriptsize $\sigma(n)$}}}
\end{picture}
\end{center}

There is a stunning example of a strongly homotopy Lie algebra
taken from Nature -- the convolution product of functionals on the
space of closed strings, 
see~\cite{saadi-zwiebach,stasheff:topology-proceedings}. 
Another place where strongly homotopy
Lie algebras naturally occur is the space of horizontal forms on the
infinite jet bundle, where they appear as lifts of brackets on the
space of Lagrangian functionals,
see~\cite{barnich-fulp-lada-stasheff:preprint,markl-shnider:eu}.
}
\end{example}

Homotopy versions of commutative associative algebras were introduced
in~\cite{kadeishvili:ttmi85} by T.~Kadeishvili. 
They are called  $\cinfty$, commutative, or
balanced $\ainfty$-algebras (compare also~\cite{markl:JPAA92}).

Let us define a map $\alpha_{\Ass} : \ainfty \to \Ass$ by
\[
\alpha_{\Ass}(\ \ \oper\ )= \mbox{  the  multiplication }\mu \in \Ass(2),
\] 
while $\alpha_{\Ass}$ is trivial on the remaining generators of 
$\ainfty$. 
We recommend as an exercise to verify that 
$\alpha_{\Ass} \circ \pa = 0$, which means that $\alpha_{\Ass}$ is a
map of differential operads, if we interpret $\Ass$ as a differential
operad with trivial differential.
It can be shown that $\alpha_{\Ass}$ induces an
isomorphism of homology, $H_*(\ainfty,\pa) \cong\Ass$. 
There is also a map $\alpha_{\Lie} :
\linfty \to \Lie$ having similar properties.

\section{Concepts}

Let us try to sum up common features of 
the above examples. In both cases, homotopy algebras were algebras
over a
differential graded operad which was free as an operad. 
Moreover, we had 
natural maps $\alpha_{\Ass} :(\ainfty,\pa) \to (\Ass,0)$ and 
$\alpha_{\Lie}: (\linfty,\pa) \to (Lie,0)$ inducing isomorphisms in  
homology. 

Recall that, in~\cite{markl:zebrulka}, 
we formulated the following definition.

\begin{definition}
Let $\P$ be an operad in the category of differential graded vector
spaces. A minimal 
model of $\P$ is a differential graded operad 
$\M_{\P}= (\F(E),\pa)$, where ${\cal F}(E)$ is the free operad on
a collection $E$, 
together with a map $\alpha_{\P} : \M_{\P} \to \P$ 
that is a homology isomorphism.
The minimality means that we assume
that $\pa(E)$ consists of decomposable elements of the free operad
$\F(E)$. 
\end{definition}

We proved, for each differential graded operad $\P$,
the existence and a kind of uniqueness of the minimal model. 
When $\P$ has trivial
differential (which was the case of all our examples), 
the minimal model is in fact bigraded, but we will not use this property.
Having in mind the above examples, we proposed,
in~\cite{markl:zebrulka}, 
the following definition.

\begin{definition}
\label{h_alg}
A homotopy $\P$-algebra is a differential graded vector space $A =
(A,\pa)$ with an action of the minimal model $\M_{\P}$ of $\P$, in other
words, a differential graded algebra over the operad $\M_{\P}$.
\end{definition}

For so called Koszul quadratic operads, the minimal model coincides
with the cobar construction on the quadratic dual of the
operad, and our definition agrees with the one given by V.~Ginzburg
and A.~Kapranov in~\cite{ginzburg-kapranov:DMJ94}.

Our approach keeps the symmetric group action (i.e.,~on the level of
algebras,
the commutativity,
anticommutativity, etc.) strict, so we
cannot handle algebraic variants of homotopy-everything operads
such as those considered 
in~\cite{smirnov:USSRSb82,smith:memoirs}, 
nor various other homotopy
algebras where the symmetry is relaxed up to a homotopy -- like homotopy
Gerstenhaber algebras in the sense of M.~Gerstenhaber and
A.A.~Voronov~\cite{gerstenhaber-voronov:FAP95}.

In the following section we try to indicate what must be done 
to justify our approach.

\section{Perspectives}
\label{lights}

There exists a natural concept of a map ${\bf
f} :A \to B$ of homotopy $\P$-algebras. It is a system of
degree zero multilinear maps ${\bf f} =\{f_n : A^{\ot n}\to B\}_{n\geq
1}$ satisfying a certain set of axioms. 
One of the axioms postulates that $f_1: (A,\pa)\to (B,\pa)$ is a
morphism of the underlying differential graded vector spaces. Let
$h\P$ denote the category of homotopy $\P$-algebras in the sense of
Definition~\ref{h_alg} and their maps.

\begin{definition}
Let $A$ and $B$ be two homotopy $\P$-algebras and $g : (A,\pa)\to
(B,\pa)$ a map of underlying differential graded vector spaces. A
homotopy $\P$-structure on $g$ is an $\hpalg$-map ${\bf
f}:A \to B$ such that $f_1 =g$.
\end{definition}

In order to justify the notion of homotopy $\P$-algebras and their
maps, 
we need
to show that homotopy $\P$-structures are `stable under a homotopy.'
The precise meaning of this was explained, for topological
spaces, in Chapter~1 of~\cite{boardman-vogt:73}.  
Corresponding statements translated to our algebraic
language are the following.
\begin{enumerate}
\item
For each homotopy $\P$-algebra $A$, differential graded vector space
$B = (B,\partial)$ and a map $g: (A,\pa) \to (B,\pa)$ that is a
homology isomorphism, there exist a homotopy $\P$-structure on
$(B,\pa)$ and a homotopy $\P$-structure ${\bf f}:A \to B$
on $g$.
\item
Suppose $A$ and $B$ are two homotopy $\P$-algebras and ${\bf f}:A \to
B$ a homotopy $\P$-algebra map. Suppose that $g:(A,\pa) \to (B,\pa)$
is a differential map that induces the same map of homology as
$f_1$. Then there exists a homotopy $\P$-structure on $g$.
\item
Suppose that ${\bf f}:A \to
B$ is a homotopy $\P$-algebra map such that $f_1$ is a homology
isomorphism. Suppose that  $g:(B,\pa) \to (A,\pa)$ is a homology
inverse of $f_1$. Then there exists a homotopy $\P$-structure on $g$.
\end{enumerate}

\noindent 
We recommend, as an easy exercise, to prove that conditions 1~and 2
imply
\begin{itemize}
\item[$1'.$]
For each homotopy $\P$-algebra $B$, differential graded vector space
$A = (A,\partial)$ and a map $g: (A,\pa) \to (B,\pa)$ that is a
homology isomorphism, there exist a homotopy $\P$-structure on
$(A,\pa)$ and a homotopy $\P$-structure ${\bf f}:A \to B$
on $g$.
\end{itemize}

As far as we know, nobody has considered the
above properties in full generality, though there are several partial
results in this direction. Let us quote at least the following
theorem due to 
T.~Kadeishvili~\cite[p.~232]{kadeishvili:RMS80}, 
which is a special case of $1'$~for the
category of $\ainfty$-algebras, with $(A,\pa) = (H(C),\pa=0)$ and $B =
(C,\pa,\mu_2,0,\ldots)$. 

\begin{theorem}
Let $(C,\partial)$ be a chain algebra. Then there exists an
$\ainfty$-structure $\{ X_k;k\geq 2\}$ on the graded space $H(C)$
having 
$X_2(a,b)=a\cdot b$ (the multiplication induced by $\mu_2$), 
together with an $\ainfty$-homomorphism ${\bf f}:
(H(C),0,X_2,X_3,\ldots)
\to(C,\pa,\mu_2,0,\dots)$ such that $f_1 : H(C) \to C$ is a
homology isomorphism.
\end{theorem}

A similar statement for balanced $\ainfty$-algebras was 
proved in~\cite{markl:JPAA92} by the 
author.
Another direction of results that indicates a certain homotopy stability of
homotopy algebraic structures is represented by various versions of
the Perturbation Lemma, 
see~\cite{gugenheim-stasheff,huebschmann-kadeishvili:MZ1991}.


\vskip3mm
\catcode`\@=11
\noindent
Mathematical Institute of the Academy, \v Zitn\'a 25, 115 67
Praha 1, Czech Republic,\hfill\break\noindent
email: {\tt
markl@math.cas.cz}\hfill\break\noindent

\end{document}